\documentclass{amsart}
\usepackage{color}
\usepackage{url}
\usepackage{amssymb,url}
\usepackage{array}
\usepackage[foot]{amsaddr}
\usepackage[norelsize,lined,boxed,commentsnumbered]{algorithm2e}

\usepackage{graphicx}
\usepackage{caption}
\usepackage{subcaption}

\graphicspath{{./data/}}

\title
[Inverse moment problem for convex polytopes: implementation]
{The inverse moment problem for convex polytopes: implementation aspects}

\author[Nick Gravin, Danny Nguyen, Dmitrii Pasechnik, Sinai Robins]
{Nick Gravin$^1$,  Danny Nguyen$^2$,
Dmitrii V. Pasechnik$^3$, Sinai Robins$^4$}
\address{$^1$ Microsoft Research, 1 Memorial Drive, Cambridge, MA 02142, USA
\newline
$^2$ Department of Mathematics, UCLA, Box 951555,
Los Angeles, CA 90095-1555, USA.
\newline
$^3$ Department of Computer Science, University of Oxford, Wolfson Building,
Parks Road, Oxford OX1 3QD, UK.
\newline
$^4$ School of Physical and Mathematical Sciences, Nanyang
Technological University, 21 Nanyang Link, 637371 Singapore.
\newline
{\rm Supported by Singapore Ministry of Education ARF Tier 2 Grant MOE2011-T2-1-090.}}
\date\today

\theoremstyle{definition}
\newtheorem{prop}{Proposition}[section]

\newtheorem{cor}{Corollary}

\newtheorem{rem}[prop]{Remark}

\theoremstyle{remark}
\newcounter{reml}

\newcommand{\C}{{\mathbb{C}}}

\newcommand{\K}{{\text{Ker}(\mathbf{H})}}

\newcommand{\Q}{{\mathbb{Q}}}
\newcommand{\R}{{\mathbb{R}}}
\newcommand{\V}{{\text{Vert}}}

\renewcommand{\epsilon}{\varepsilon}
\renewcommand{\theta}{\vartheta}

\newcommand{\x}{\mathbf{x}}

\newcommand{\z}{\mathbf{z}}

\newcommand{\bv}{\mathbf{v}}

\newcommand{\bH}{\mathbf{H}}

\newcommand{\bil}[2]{\langle{#1},{#2}\rangle}

\newcommand{\marginnote}[1]{\vrule width0pt height0pt depth0pt
 \vadjust{\vbox to0pt{\vss\hbox to\hsize{\hskip\hsize\quad
 #1\hss}\vskip1.5pt}}}

\newcommand{\beq}{\begin{equation}}
\newcommand{\eeq}{\end{equation}}

\newcommand{\BrLa}{BBaKLP}

\begin{document}
\begin{abstract}
We give a detailed technical report on the implementation of the
algorithm presented in \cite{GLPR12} for reconstructing an $N$-vertex
convex polytope $P$ in $\R^d$ from the knowledge of $O(Nd)$ its moments.

\end{abstract}

\maketitle

\section{Problem description}

Our main object of interest is a convex polytope $P\subset\R^d$ with $N$-vertices.
We assume that the polytope $P$ has a polynomial density $\rho(\x)$ defined in
the interior of $P$. For any multivariate polynomial $g(\x)$ the corresponding
moment $\mu_g$ of $P$ is given by

\[
\mu_g:= \int_P  g(\x)\cdot \rho(\x) d\x.
\]

We note that if all vertices of $P$ are rational (have rational coordinates) and
$\rho\in\Q[\x]$, then every moment $\mu_g$ of $P$
for a polynomial $g\in\Q[\x]$ is a rational number as well. Why this is true will
become clear in the next section.

\paragraph{\bf Input}  As an input to our problem we receive $O(Nd)$ moments
of some underlying $N$-vertex convex polytope $P\subset \R^d$.
\paragraph{\bf Output} The goal is to reconstruct $P$ (coordinates of the vertices).

In our computational experiments we did a few simplification assumptions about the underlying polytope:
\begin{enumerate}
\item we work with uniform density, i.e., $\rho(\x)=1$ for any $\x\in P$;
\item we focus on {\em simple} polytopes, i.e., polytopes where each vertex has exactly $d+1$ incident edges.
\end{enumerate}

The latter assumption is equivalent to saying that $P$ is a generic polytope in a hyper-plane description of the polytope, i.e.,
no $d+1$ supporting hyperplanes of $P$ have common intersection. In order to construct
a random simple polytope our computational experiments we intersect a few half spaces each supported by a randomly chosen hyperplane.

We considered the problem in two different models of arithmetic:
\begin{enumerate}
\item vertices of $P$ are rational and rational moments are given in the input exactly;
\item vertices of $P$ have real coordinates and moments are given with certain precision.
\end{enumerate}

\section{Preliminaries}

For a non-negative integer $j$ the $j$-th {\em axial moment}
of $P$ in the direction $\z\in\R^d$ with respect to density $\rho$ is given by
\[
\mu_j(\z):= \mu_{j,\rho} (\z) := \int_P  \langle \x, \z \rangle  ^j \rho(\x) d\x.
 \]

We remark that $\langle \x, \z \rangle^j$ is a homogeneous polynomial of degree $j$
for any fixed direction $\z$.

Let the set of all vertices of $P$ be given by $\V(P)$.
For each $\bv \in \V(P)$, we consider a fixed set of
vectors, parallel to the edges of $P$ that are incident with $\bv$, and call
these edge vectors $w_1(\bv)$,\dots $w_d(\bv)$.  Geometrically, the polyhedral
cone generated by the non-negative real span of these edges at $\bv$ is called
the tangent cone at $\bv$, and is written as $K_\bv$.  For each simple tangent
cone $K_\bv$, we let $|\det K_\bv|$ be the volume of the parallelepiped formed
by the $d$ edge vectors $w_1(\bv), \dots, w_d(\bv)$.  Thus, $|\det K_\bv| = |
\det( w_1(\bv),  \dots, w_d(\bv)) |$, the determinant of this parallelepiped.

The following results of \BrLa\ \cite{MR1079024} tells us
\beq\label{brion} \mu_j (\z)= \frac{j! (-1)^d}{ (j+d)!}
\sum_{\bv\in \V(P)} \bil{\bv}{\z}^{j+d} D_\bv(\z), \text{ where}
\eeq
\beq\label{Dvz} D_\bv(\z):=\frac{|\det
K_\bv|}{\prod_{k=1}^d\bil{w_k(\bv)}{\z}},
\eeq
for each $\z \in \R^d$ such that the denominators in $D_\bv(\z)$ do not vanish.
Moreover,
\beq\label{brionzero} 0=\sum_{v\in \V(P)} \bil{\bv}{\z}^{j}
D_\bv(\z), \text{ for each } 0 \leq j \leq d-1.\eeq

In particular, from \eqref{brion},\eqref{Dvz} it is easy to see that every moment $\mu_j(\z)$
is a rational number, if $P$ is a rational polytope and $\z\in\Q^d$. Since any polynomial $g\in\Q[\x]$
can be expressed as a rational linear combination of the powers of linear forms with rational coefficients, we
can conclude that $\mu_g\in\Q$.

Rewriting the above equations in the matrix form we get

\beq\label{momentmatrix}
\begin{pmatrix}
1&1&\dots&1\\
 \langle \bv_1, \z \rangle &  \langle \bv_2, \z \rangle & \dots &  \langle \bv_N, \z \rangle   \\
{ \langle \bv_1, \z \rangle}^2 & { \langle \bv_2, \z \rangle}^2 & \dots &  {\langle \bv_N, \z \rangle}^2  \\
\vdots&\vdots&\dots&\vdots \\
  { \langle \bv_1, \z \rangle}^k & { \langle \bv_2, \z \rangle}^k & \dots &  {\langle \bv_N, \z \rangle}^k  \\
\end{pmatrix}
\begin{pmatrix}
D_{\bv_1}(\z)\\
\vdots\\
D_{\bv_N}(\z)
\end{pmatrix}=
\begin{pmatrix}
c_0\\
\vdots\\
c_{k}
\end{pmatrix},
\eeq
where
\beq\label{c-vector}
\left (c_0, \dots, c_{k} \right) = \left(0, \dots, 0, \frac{d!(-1)^d}{0!} \mu_0,  \frac{(1+d)!(-1)^d}{1!} \mu_1, \dots, \frac{k!(-1)^d}{(k-d)!} \mu_{k-d}\right),
\eeq
so that the vector ${\bf{c}} = \left(c_0, \dots, c_{k}\right)$ has zeros in the first $d$ coordinates, and scaled moments in the last $k+1-d$ coordinates.

For a fixed $m\ge N+1$ let:

\beq
\bH (c_0,\dots,c_{2m-2}):=
\begin{pmatrix}
c_0&c_1&\dots&c_{m-1}\\
c_1&c_2&\dots&c_{m}\\
\vdots&\vdots&\dots&\vdots\\
c_{m-}&c_{m+1}&\dots&c_{2m-2}
\end{pmatrix}.
\eeq

Below is given the algorithm (a variant of the Prony method) from \cite{GLPR12} of how to find the projections
of vertices of $P$ onto a general position axis $\z\in \R^d$.

\begin{algorithm}[H]
\begin{enumerate}

\item  Given $2m-1 \ge 2N+1$ moments $c_0,\dots,c_{2m-2}$ for $\z$, construct \\
\noindent a square Hankel matrix $\bH(c_0,\dots,c_{2m-2}).$
\medskip
\item Find the vector $v=\left(a_0, \ldots, a_{M-1}, 1, 0, \ldots, 0 \right)$ in $\K$ \\
\noindent with the minimal possible $M.$ It turns out that the number of \\
       \noindent vertices $N=M$.
\medskip
\item The set of roots $\{x_i(\z) = \bil{\bv_i}{\z} |\bv_i\in\V(P)\}$ of polynomial
      $p_\z(t) = a_0 + a_1t + \ldots + a_{N-1} t^{N-1} + t^N$
       then equals the set of \\
       \noindent projections of $\V(P)$ onto $\z$.
\end{enumerate}
\caption{Computing projections.}\label{fig:compproj}
\end{algorithm}

Next the algorithm in \cite{GLPR12} finds projections of the vertices on $d$ different linearly independent directions $\z\in\R^d$ and
matches the projections on the first direction with the projections on each of the rest $d-1$ directions. In order to do each matching
between the first $\z_1$ and $i$-th $\z_i$ directions, vertex projections on a new direction $\z_{1i}$ in the plane
spanned by $\z_1$ and $\z_i$ are reconstructed. These extra projections on the direction $\z_{1i}$ allow to restore the right matching between
the projections on $\z_1$ and $\z_i$ with very high probability.

\section{Actual Implementation}
Our implementation was done in Sage \cite{sage}.
\paragraph{\bf Reconstructing projections on $\z$.}Coming to the main part, we deviated a little bit from our original Prony
method in computing the axial projections. Namely,  we do not go directly
on finding the kernel of the Hankel system but look at the problem from the
perspective of Pade approximation instead. The moments can be viewed as
coefficients in the expansion of a rational function, which we can approximate
if enough data is known. Specifically, recalling \eqref{momentmatrix} and \eqref{Dvz},
we may write the following univariate generating function for the sequence of scaled moments $\{c_k\}$

\begin{eqnarray}
\sum_{k=0}^{\infty}c_{k}t^k &=& \sum_{k=0}^{\infty}t^k\sum_{i=1}^{N} \langle\bv_i,\z\rangle^k D_{\bv_i}(\z)\notag\\
&=&\sum_{i=1}^{N}D_{\bv_i}(\z) \sum_{k=0}^{\infty}t^k\langle\bv_i,\z\rangle^k=
\sum_{i=1}^{N}\frac{D_{\bv_i}(\z)}{1-t\langle\bv_i,\z\rangle}.
\label{pade}
\end{eqnarray}

Therefore, $c_k$ are the coefficients in the Taylor series expansion of
$p_\z(t)/q_\z(t)$, where $q_\z(t) = \prod\limits_{\bv\in\V(P)} (1-t\langle\bv,\z\rangle)$ and $q_\z(t)$ is a polynomial of degree at most $N-1$.
If enough moments are known for a fixed direction $\z$ (in our case $2N$ are sufficient) then $p_\z$ and $q_\z$ can be computed. Then the roots of
$q_\z$ will give us the desired projections.

In our implementation we used one of the Pade approximation methods implemented in Sage. This is basically \texttt{scipy.misc.pade} with control of the measured moments' precision. If:
$$
\frac{p(t)}{q(t)} = \frac{a_0 + a_1 t + \dots + a_\ell t^\ell}{b_0 + b_1 t + \dots + b_m t^m} = c_0 + c_1 t + \dots + c_n t^n + \dots
$$
where $n = \ell+m$, $q_0=1$ and $c_0,\dots,c_n$ are moments then we do the following:

\begin{algorithm}[H]
\begin{enumerate}

\item  Trim the data $(c_0,\dots,c_n)$ to $k$-bit precision with $k$ specified. \\
\medskip
\item Create a matrix $C_{m\times m}$ with $C_{ij}=c_{\ell+i-j}$.
\medskip
\item Solve the system $C\cdot x = y$ with $x = (b_1,\dots,b_m)^T$ and \\
\noindent $y = -(c_{\ell+1},\dots,c_{\ell+m})^T$.
\end{enumerate}
\caption{Pade approximation.}\label{fig:pade}
\end{algorithm}

\paragraph{\bf Matching projections on different directions.}
We implemented a different and much more reliable matching procedure than the one
described in the original paper. Below we give a detailed description of the new matching method.

As was remarked in \cite{GLPR12}, formulas \eqref{brion}
and \eqref{Dvz} are valid not only for $\z\in\R^d$ but also for $\z\in\C^d$. The latter
means that every point in $P\subset\R^d$ and each $\bv\in\R^d$ are regarded as
complex vectors with all zero imaginary components and $\bil{\bv}{\z}$ is regarded as
a standard sesquilinear inner product in $\C^d$.

Thus, we also can write \eqref{momentmatrix} for complex $\z=\z_{re}+i\cdot\z_{im},$
where $\z_{re},\z_{im}\in\R^d$.
We observe that $\bil{\bv}{\z}=\bil{\bv}{\z_{re}}+i\cdot\bil{\bv}{\z_{im}}$ and
\begin{eqnarray}
\mu_j(\z) &=&  \int_P  \Big(\langle \x, \z_{re} \rangle + i\cdot
\langle \x, \z_{im} \rangle \Big)^j \rho(\x) d\x  \notag\\
&=&\int_P  g_1\Big(\langle \x, \z_{re} \rangle,
\langle \x, \z_{im} \rangle \Big) \rho(\x) d\x+ i\cdot\int_P g_2\Big(\langle \x, \z_{re} \rangle,
\langle \x, \z_{im} \rangle \Big)\rho(\x) d\x, \notag
\end{eqnarray}
where $g_1$ and $g_2$ are homogeneous real polynomials in two variables of degree $j$.
Hence, by receiving in the input moments $\mu_{g_1},$ $\mu_{g_2}$ we may find $\mu_j(\z)$ for any $\z\in\C^d$.

We further may write \eqref{pade} for $\z\in\C^d$ and find $q_\z(t)$ from the first $2N$ moments. Next we
find all complex roots of the polynomial $q_\z(t)$, which give us already matched projections on $\z_{re}$ and $\z_{im}$.
In our algorithm we fix some general position vector $\z_{re}\in\R^d$ and consider $d-1$ vectors $\z_j\in\R^d$, such that $\z_{re}$ and all
$\z_j$ are linearly independent. We match projections on $\z_{re}$ with the projections on $\z_j$ by taking $\z_{im}=\z_j$ for each $j$.
A big advantage of this matching method is that it is much less prone to numerical errors. In particular, for $d=2$ this method will provide
an answer in any case, in other word for $d=2$ our problem is well posed. For $d\ge 3$ there might be a problem that projections on $\z_{re}$
are different when we match them with projections on different $\z_j$. We simply get around this problem by using an ascending order over projections on $\z_{re}$ each time when we do such a matching.

\begin{rem}Interestingly, if we fix unit and orthogonal to each other directions $\z_{re}$ and $\z_{im}$,
then polynomials $g_1\Big(\langle \x, \z_{re} \rangle,
\langle \x, \z_{im} \rangle \Big)$ and $g_2\Big(\langle \x, \z_{re} \rangle,
\langle \x, \z_{im} \rangle \Big)$ considered as multivariate polynomials of $\x$ are harmonic functions, i.e., $\Delta g_1(\x)=\Delta g_2(\x)=0.$
One can read more on harmonic moments in e.g. \cite{PS12}.
\end{rem}
\begin{proof} We recall that Laplace operator $\Delta$ is invariant under the isometry group of $\R^d$. Therefore, we may assume
that $\z_{re}$ is simply the first coordinate vector of $\x$ and $\z_{im}$ is the second coordinate vector of $\x$. Now we need only to verify that
$\Delta=\frac{\partial^2}{\partial x^2}+\frac{\partial^2}{\partial y^2}$ when applied to the real and imaginary part of $(x+i\cdot y)^j$ is zero.
Indeed, we have
\[
\Delta(x+i\cdot y)^j=j(j-1)(x+i\cdot y)^{j-2}+j(j-1)i\cdot i\cdot(x+i\cdot y)^{j-2}=0.
\]
\end{proof}

\begin{cor}From harmonic moments only, one may reconstruct vertices of a convex polytope $P$.
\end{cor}

\section{Numerical Experiments}

We did our numerical experiments first in the exact arithmetic, i.e with rational precision, to test
the exact algorithm from \cite{GLPR12} and adjust the part of the algorithm
for selecting random directions. In this mode our implementation was far
from optimal in terms of running time with exact arithmetic. The reason
for that is due to the inherent limitation of Sage's rational arithmetic.

\begin{table}[h]
\begin{tabular}{|c|c|c|c|c|}
\hline
Dimension & Number of & Exact      & Float      & Allowed\\
          & vertices  & Arithmetic & Arithmetic & Error     \\ \hline
2         & 10        & 0.47 sec   & 0.07 sec   & E-3           \\ \hline					
3         & 20        & 39 sec     & 0.42 sec   & E-3           \\ \hline					
4         & 30        & $>$ 5 mins   & 1.89 sec   & E-3           \\ \hline				
5         & 40        & $>$ 10 mins   & 7.10 sec  & E-3           \\ \hline 					
\end{tabular}
\caption{Efficiency Benchmarks.}
\label{table:timing}
\end{table}

\begin{figure}[h]
\centering
    \begin{subfigure}[b]{0.4\textwidth}
        \centering
        \includegraphics[width=\textwidth]{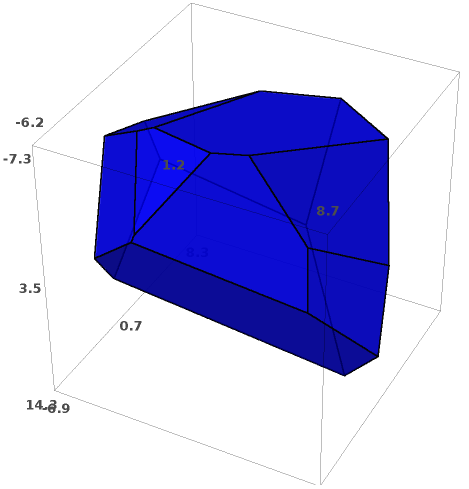}
        \caption{Original Polyhedron}
        \label{fig:3dim_20_rational}
    \end{subfigure}%
    \qquad
    \begin{subfigure}[b]{0.4\textwidth}
        \centering
        \includegraphics[width=\textwidth]{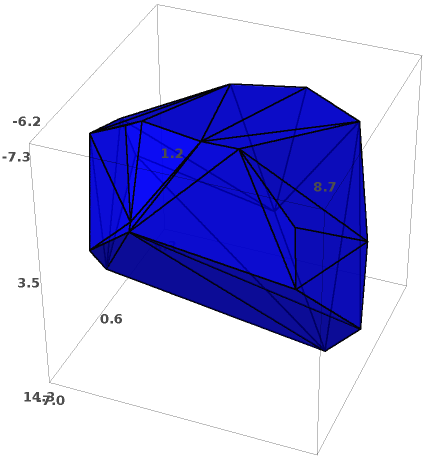}
        \caption{75-digit Float Arithmetic}
        \label{fig:3dim_20_float75}
    \end{subfigure}
\caption{3D polyhedron with 20 vertices}
\label{fig:3dim_20}
\end{figure}

However, as can be seen from the table \ref{table:timing}, converting numerical data into float precision yields drastic improvements in terms of running time. The allowed error on recovered projections is small enough and leaves the shape almost intact. On the figure \ref{fig:3dim_20} are two images of the same 20-vertex polyhedron reconstructed with rational and float arithmetic. Note that with float arithmetic, tiny errors in projections altered co-planarity of many vertices and thus many facets are triangulated although the general shape is still preserved.

\begin{table}[h]
\begin{tabular}{|c|c|c|c|}
\hline
Number of & Error of  & Error or   & Error of \\
vertices  & order E-3 & order E-6  & order E-9 \\ \hline
4         & 20 bits   & 25 bits    & 35 bits \\ \hline					
8         & 30 bits   & 40 bits    & 45 bits \\ \hline					
12        & 45 bits   & 55 bits    & 65 bits \\ \hline				
16        & 60 bits   & 65 bits    & 75 bits \\ \hline 					
20        & 75 bits   & 80 bits	   & 90 bits \\ \hline
40        & 160 bits  & 170 bits   & 210 bits \\ \hline
\end{tabular}
\caption{Errors v.s. Float Precision}
\label{table:trade-off}
\end{table}

When noise is introduced to the measured moments, exact arithmetic becomes inapplicable. Float arithmetic on the other hand can tolerate errors to some degree. However, to retrieve projections with high precision, our method turns out to be very sensitive. In table \ref{table:trade-off} we compare the precision level required with float arithmetic versus error tolerability.

\begin{figure}[h]
\centering
    \begin{subfigure}[b]{0.4\textwidth}
        \centering
        \includegraphics[width=\textwidth]{3D_20vert_rational}
        \caption{Original Polyhedron}
    \end{subfigure}%
    \qquad
    \begin{subfigure}[b]{0.4\textwidth}
        \centering
        \includegraphics[width=\textwidth]{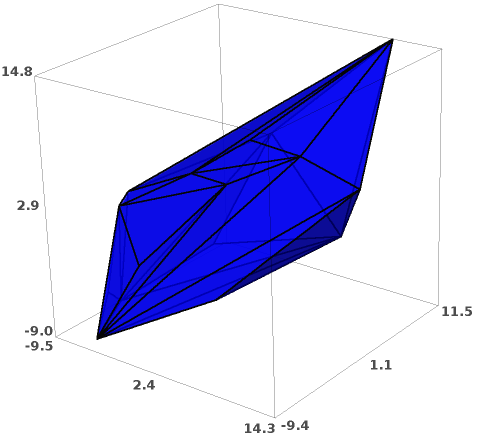}
        \caption{Distorted Recovery}
        \label{fig:3dim_20_float60}
    \end{subfigure}
\caption{Float arithmetic: 60 digits precision}
\label{fig:distorted_recovery}
\end{figure}

We give an example of insufficient precision that results in distortions of the reconstructed shape. With the previous 20-vertex polyhedron where moments are measured now to only 60 bits of precision, the recovered shape looks as is shown on the figure \ref{fig:distorted_recovery}.

We would like to remark that the use of complex moments improved precision a lot compared to the real moments. Here is a concrete example of a 8-vertex polyhedron with the matrix $V$ containing vertex coordinates and $A$ representing its adjacency matrix.

\[
V=\bordermatrix{
~ &v_1    &v_2        &v_3      &v_4      &v_5      &v_6        &
v_7    &v_8    \cr
x & 17/4  & 249/121   & -719/74 & -66/43  & -82/91  & -1588/133 &
545/37 & 69/7  \cr
y & -14/3 & -211/121  & -373/74 & -267/43 & -219/91 & 414/133   &
765/37 & 59/21 \cr
z & -7/12 &  1963/121 & 426/37  & -108/43 & -148/13 & -46/133   &
-85/37 & -41/3 \cr
}
\]

\[
\text{Adjacency matrix: } \begin{pmatrix}
0 & 1 & 0 & 1 & 0 & 0 & 0 & 1\\
1 & 0 & 1 & 0 & 0 & 0 & 1 & 0\\
0 & 1 & 0 & 1 & 0 & 1 & 0 & 0\\
1 & 0 & 1 & 0 & 1 & 0 & 0 & 0\\
0 & 0 & 0 & 1 & 0 & 1 & 0 & 1\\
0 & 0 & 1 & 0 & 1 & 0 & 1 & 0\\
0 & 1 & 0 & 0 & 0 & 1 & 0 & 1\\
1 & 0 & 0 & 0 & 1 & 0 & 1 & 0
\end{pmatrix}
\]



\begin{figure}[h]
\centering
  \includegraphics[width=0.5\textwidth]{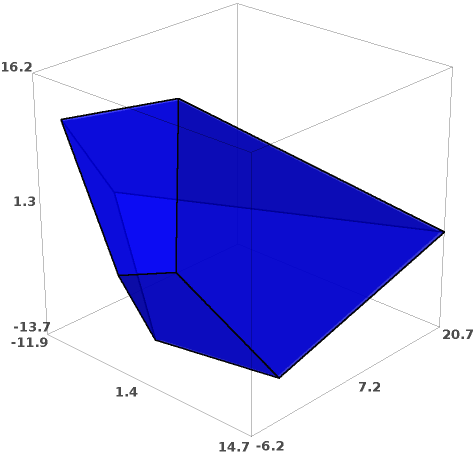}
  \caption{3D Polyhedron with 8 vertices}
  \label{fig:3dim_8vert}
\end{figure}


\texttt{z = vector([2,3,4])} is the random vectors upon which vertices are projected. The exact projections are:

\texttt{Proj:~[-54.56, -31.74, -26.52, -15.92, -7.83, 11.50, 63.78, 82.30]}


With moments measured in the real field with 25-bit precision wrapping, we recovered the projections as:

\texttt{RealField(25):~[-54.56, -30.87, -23.14, 11.46, 63.78, 82.30]}


Notice that some projections are missed, because the computations done in the Real filed with $25$-digit precision have affected slightly the coefficients of $p_\z(t)$ and, therefore, some real roots of $p_\z(t)$ have disappeared. Now with a randomly chosen complex component, we can take $\z =$\texttt{ vector([2,3,4]) + I*vector([-5,2,-8])} and carry out the same computations in the complex field with 25-digit precision and the result recovers all 8 projections with much better precision \texttt{ComplexField(25): [-54.56, -31.81, -26.48, -15.93, -7.84, 11.49, 63.78, 82.30]}


\section{Conclusions}

In our computational experiments with exact arithmetic and precise measurements we achieved the expected performance and precision guarantees and have improved the original algorithm suggested in \cite{GLPR12} in certain respects. Namely, we implemented significantly more robust matching procedure of the vertex projections by recovering projections of the vertices on a complex plane instead of a single direction recovery as was proposed in the original work; we implemented an easier and more practical procedure based on Pade approximation to recover the projections on the given complex plane and/or single real axis. One of the interesting implications of the former methodology is that harmonic moments (polynomials $p(\x)$, s.t. $\Delta p(x) = 0$) are sufficient to recover vertices of any convex polytope as well as vertices of a non-convex polytope, if the respective coefficients at the vertices do not vanish. We note that there are examples of different non convex polytopes with exactly the same set of harmonic moments.

On the negative side, in the numerical experiments with bounded precision we have seen a very high sensitivity of our methodology to numerical inaccuracies.
The latter is an unavoidable obstacle to the practical usage of our algorithm.


\bibliographystyle{alpha}
\bibliography{../poly,../master}

\end{document}